\documentclass[12pt]{article}
\usepackage{epsfig}
\usepackage{amsfonts}
\usepackage{amssymb}
\usepackage{amsmath}
\usepackage{amsthm}
\usepackage{latexsym}
\usepackage{graphicx}
\usepackage{bm}
\usepackage{tabularx}
\usepackage{booktabs} 
\usepackage{enumerate}
\usepackage{caption}
\usepackage[titletoc]{appendix}
\usepackage[numbers]{natbib}
\usepackage{bbm}
\usepackage{url}
\usepackage{hyperref}
\usepackage{color}
\usepackage{xcolor}
\usepackage[section]{placeins}
\usepackage{subcaption}
\usepackage{subfloat}
\usepackage{placeins}

\usepackage{floatrow}
\newfloatcommand{capbtabbox}{table}[][\FBwidth]

\usepackage{blindtext}

\DeclareUnicodeCharacter{0308}{ü}

\let\Oldsection\section
\renewcommand{\section}{\FloatBarrier\Oldsection}

\let\Oldsubsection\subsection
\renewcommand{\subsection}{\FloatBarrier\Oldsubsection}

\let\Oldsubsubsection\subsubsection
\renewcommand{\subsubsection}{\FloatBarrier\Oldsubsubsection}

\catcode`\@=11

\theoremstyle{plain}

\theoremstyle{definition}

\theoremstyle{remark}

\begin{document}
\title{A Note on the Interpretation of Distributed Delay Equations}
\author{ Philip Doldo \\ Center for Applied Mathematics \\ Cornell University
\\ 657 Rhodes Hall, Ithaca, NY 14853 \\  pmd93@cornell.edu  \\ 
\and
  Jamol Pender \\ School of Operations Research and Information Engineering \\ Center for Applied Mathematics \\ Cornell University
\\ 228 Rhodes Hall, Ithaca, NY 14853 \\  jjp274@cornell.edu 
 }    

\maketitle

\begin{abstract}

Distributed delay equations have been used to model situations in which there is some sort of delay whose duration is uncertain. However, the interpretation of a distributed delay equation is actually very different from that of a delay differential equation with a random delay. This work explicitly highlights this distinction as it is an important consideration to make when modeling delayed systems in which the delay can take on several values.

\end{abstract}


\section{Introduction} \label{sec_intro}

Delay differential equations (DDEs) are often used to model systems in which some sort of time delay is present. While much of the literature on DDEs is concerned with models that have constant delays, a significant portion of the literature considers models in which the delay can vary randomly as it is realistic in many applications that the delay present in a system may not always be the same. One way in which this randomness has been modeled in the literature is through distributed delay equations, which can roughly be viewed as DDEs where the delayed terms in the model form an expectation. For example, the DDE that we may genuinely be interested in understanding could be 
\begin{eqnarray}
\overset{\bullet}{u}(t) &=& \alpha \cdot u(t - \Delta) \label{eqn_random}
\end{eqnarray}

\noindent where $\Delta$ is some nonnegative random variable and the corresponding distributed delay equation which would be analyzed is 
\begin{eqnarray}
\overset{\bullet}{u}(t) &=& \alpha \cdot \mathbb{E}[ u(t - \Delta)]
\end{eqnarray}

\noindent where the expectation is over the distribution of the random variable $\Delta$ and we note that this reduces to a multi-delay differential equation when $\Delta$ has a discrete probability distribution. Examples of work that has analyzed distributed delay equations include \cite{campbell2009approximating, tang2004oscillation, pender2017queues, novitzky2019nonlinear, bernard2001sufficient, moruarescu2007stability, cassidy2020distributed, doldo2020multi, novitzky2020queues, jessop2010approximating, novitzky2020limiting, berezansky2001oscillation, braverman2012absolute, cao2019bifurcations}. However, many of these works use distributed delay equations in a way that is misleading. Consider the following quote from \cite{jessop2010approximating} which references a distributed delay equation model: "But many biological and physical events, such as regeneration, recovery period from a disease, or signal conduction, may not take exactly the same time in each instance. Hence model (2) can be further improved by including a distribution of delays".  This quotation from \cite{jessop2010approximating} makes it seem that the delay that they want to model, such as a regeneration or recovery period, is not a constant event and is in some sense random.  Thus, they propose to use distributed delay models to capture the impact of this randomness.  We will show below that ignoring how we capture this randomness can lead to incorrect dynamical systems models.  

When analyzing some random process such as the $u(t)$ in Equation \ref{eqn_random}, it is often useful to try to understand the mean dynamics of the system which can experimentally be approximated by sampling many realizations of the random process and averaging. We believe that it is important to make it clear that the information obtained from sampling is \emph{not} the same as the information obtained by analyzing a distributed delay equation as an explicit mention of this distinction is absent from the literature and should be considered when justifying which models to use to describe systems with random delays. In this brief paper we focus on clarifying this distinction and discussing the interpretations of sampling and distributed delay equations.  We hope to elucidate when it is appropriate to use a distributed delay in equation in modeling random phenomena.  



\subsection{Sampling DDEs with a Random Delay}

In this section, we briefly discuss what it means to sample a DDE that has a delay which is a random variable. That is, we will discuss what it means to sample solutions $v(t)$ to the DDE
\begin{eqnarray}
\overset{\bullet}{v}(t) &=&  \alpha \cdot v(t - \Delta) \label{random}
\end{eqnarray}

\noindent where $\Delta$ is a nonnegative random variable. We can view the solution $v(t)$ as a random variable for which a given realization of depends on a realization of $\Delta$. Thus, to sample realizations of $v(t)$, we will sample realizations of $\Delta$ and solve DDEs corresponding to each sample value of $\Delta$. 

Suppose we sample $M$ realizations of the random variable $\Delta$ and for the $i^{\text{th}}$ sample it takes on the value $\Delta = \Delta^{(i)}$. Corresponding to this sample is the DDE 
\begin{eqnarray}
\overset{\bullet}{v}^{(i)}(t) &=& \alpha \cdot v^{(i)}(t - \Delta^{(i)})
\end{eqnarray}

\noindent which can be solved to obtain $v^{(i)}(t)$, for $i  = 1, ..., M$. We denote the sample mean of these solutions when using $M$ samples by $v_{R_{M}}(t)$ which is defined as 
\begin{eqnarray}
v_{R_{M}}(t) &:=& \frac{1}{M} \sum_{i=1}^{M} v^{(i)}(t)
\end{eqnarray}

\noindent and we also define 
\begin{eqnarray}
v_R(t) &:=& \lim_{M \to \infty} v_{R_{M}}(t). 
\end{eqnarray}

\noindent By the strong law of large numbers (of course, noting that our samples are i.i.d.), we know that $$v_{R_M}(t) \to \mathbb{E}[v(t)] \hspace{2mm} \text{a.s.  as} \hspace{2mm} M \to \infty$$ where $v(t)$ is the solution to $$\overset{\bullet}{v}(t) = \alpha \cdot v(t - \Delta)$$ and the expected value is over the distribution of the random variable $\Delta$. We thus note that $v_R(t) = \mathbb{E}[v(t)]$ almost surely. If we define $v_{\delta}(t)$ to be the solution to $$\overset{\bullet}{v}_{\delta}(t) = \alpha \cdot v_{\delta}(t - \delta) $$ for some constant $\delta \geq 0$, then we see that (almost surely) $$v_R(t) = \mathbb{E}[v(t)] = \int_{0}^{\infty} v_{\delta}(t) f_{\Delta}(\delta) d \delta $$ where we are assuming that $f_{\Delta}$ is the probability density function (which, for simplicity, we are assuming exists) corresponding to the nonnegative random variable $\Delta$. 

Ultimately, it is clear that understanding the sample mean of solutions is a valid approach to understanding the mean dynamics of the solution to Equation \ref{random}.



\subsection{Discrete Random Delay and Distributed Delay Equations}

In this section we will explicitly showcase the distinction between sampling DDEs and distributed delay equations in a discrete setting. Let $\Delta$ be a discrete random variable so that $\Delta = \Delta_i$ with probability $p_i$ for $i=1,...,m$. Consider the following (deterministic) delay differential equations.
\begin{align}
\overset{\bullet}{v}_1(t) &= \alpha \cdot v_1(t - \Delta_1) \label{v1dde}\\
\overset{\bullet}{v}_2(t) &= \alpha \cdot v_2(t - \Delta_2)\\
&\vdots \nonumber \\
\overset{\bullet}{v}_m(t) &= \alpha \cdot v_m(t - \Delta_m) \label{vmdde}
\end{align}

\noindent Suppose we do $M$ samples of $\Delta$ where the $i^{\text{th}}$ sample corresponds to the solution to $$\overset{\bullet}{v}^{(i)}(t) = \alpha \cdot v^{(i)}(t - \Delta^{(i)}), \hspace{5mm} i=1, ..., M.$$ Each of the $M$ samples will correspond to one of the $m$ DDEs listed in Equations \ref{v1dde}-\ref{vmdde}. Let $M_i$ be the number of samples corresponding to $v_i(t)$ for $i=1, ..., m$ so that $M_1 + \cdots + M_m = M$. We then have that $$v_{R_{M}}(t) = \frac{M_1}{M} v_1(t) + \frac{M_2}{M} v_2(t) + \cdots + \frac{M_m}{M} v_m(t) $$ and we can deduce from the strong law of large numbers that $$\frac{M_i}{M} \to p_i \hspace{2mm} \text{a.s. as} \hspace{2mm} M \to \infty$$ and that $$v_R(t) = p_1 v_1(t) + p_2 v_2(t) + \cdots + p_m v_m(t)$$ almost surely. Differentiating, we see that 
\begin{align}
\overset{\bullet}{v}_R(t) &= p_1 \overset{\bullet}{v}_1(t) + p_2 \overset{\bullet}{v}_2(t) + \cdots + p_m \overset{\bullet}{v}_m(t)\\
&= \alpha \cdot \left[ p_1 v_1(t - \Delta_1) + p_2 v_2(t - \Delta_2) + \cdots + p_m v_m(t - \Delta_m)  \right].
\end{align}

\noindent We therefore see that this is the DDE whose solution, $v_R(t)$, corresponds to sampling solutions to $$\overset{\bullet}{v}(t) = \alpha \cdot v(t - \Delta)$$ for specific realizations of the random variable $\Delta$ and then averaging together the samples where we know that a proportion $p_i$ of the samples will correspond to the solution to $$\overset{\bullet}{v}_i(t) = \alpha \cdot v_i(t - \Delta_i)$$ for $i=1,..., m$ (almost surely in the limit as the number of samples goes to infininty).

Alternatively, consider the distributed delay differential equation
\begin{align}
\overset{\bullet}{v}_D(t) &= \alpha \cdot \mathbb{E}[v_D(t - \Delta)]\\
&= \alpha \cdot [p_1 v_D(t - \Delta_1) + p_2 v_D(t - \Delta_2) + \cdots + p_m v_D(t - \Delta_m)].
\end{align}

\noindent Compare the following two delay differential equations
\begin{align}
\overset{\bullet}{v}_R(t) &= \alpha \cdot \left[ p_1 v_1(t - \Delta_1) + p_2 v_2(t - \Delta_2) + \cdots + p_m v_m(t - \Delta_m)  \right]\\
\overset{\bullet}{v}_D(t) &=\alpha \cdot [p_1 v_D(t - \Delta_1) + p_2 v_D(t - \Delta_2) + \cdots + p_m v_D(t - \Delta_m)]
\end{align}

\noindent and take care to note the difference between the right-hand sides (i.e., $v_1 \neq v_2 \neq \cdots \neq v_R$ in general, so the equations are in different forms). In some sense, we can view the distributed delay differential equation as being obtained by averaging together the "models" or "right-hand-side operators" of the equations $$\overset{\bullet}{v}_i(t) = \alpha \cdot v_i(t - \Delta_i)$$ for $i = 1, ..., m$. That is, if we define the operator $L_i$ by $$L_i[w] = \alpha \cdot w(t - \Delta_i)$$ so that $$\overset{\bullet}{v}_i(t) = L_i[v_i]$$ for $i=1,...,m$, then we can rewrite the distributed delay differential equation as 
\begin{align}
\overset{\bullet}{v}_D(t) &= p_1 L_1[v_D] + p_2 L_2[v_D] + \cdots + p_m L_m[v_D]\\
&= \left(\sum_{i=1}^{m} p_i L_i \right)[v_D].
\end{align}

\noindent Thus, the distributed delay differential equation is essentially obtained by sampling over the operators $L_i$ (where a proportion $p_i$ of the samples will almost surely correspond to the operator $L_i$ in the limit as the number of samples goes to infinity) and averaging them together to get the "model" or "right-hand-side operator" for the distributed delay differential equation.

This highlights the important distinction between sampling \emph{solutions} of DDEs (which corresponds to $v_R(t)$) and sampling \emph{models} of DDEs (which corresponds to $v_D(t)$). In particular, sampling \emph{solutions} makes sense when there is uncertainty in the delay and sampling \emph{models} makes sense when there is uncertainty in the model.

\section{Conclusion}
\label{conclusion}

Many systems which have delays that have some uncertainty can be modeled by DDEs with random delays. In this paper we explicitly pointed out the distinction between analyzing the average behavior of a DDE with a random delay and analyzing a distributed delay equation. In particular, the former can be viewed as averaging DDE solutions whereas the latter can be viewed as averaging DDE operators. It is important to be aware of this distinction when deciding upon how to model randomness in delayed systems.  Thus, with the rise in research in uncertainty quantification, distributed delay equations should be thought of in that context as they are uncertainty with respect to the dynamical system model.  We hope that this note clarifies how people interpret distributed delay equations in the future and helps the reader understand that they should not viewed as a way of randomizing the delay.

\bibliographystyle{plainnat}
\bibliography{distributed_delay}

\begin{thebibliography}{14}
\providecommand{\natexlab}[1]{#1}
\providecommand{\url}[1]{\texttt{#1}}
\expandafter\ifx\csname urlstyle\endcsname\relax
  \providecommand{\doi}[1]{doi: #1}\else
  \providecommand{\doi}{doi: \begingroup \urlstyle{rm}\Url}\fi

\bibitem[Berezansky and Braverman(2001)]{berezansky2001oscillation}
Leonid Berezansky and Elena Braverman.
\newblock On oscillation of equations with distributed delay.
\newblock \emph{Zeitschrift f{\"u}r Analysis und ihre Anwendungen}, 20\penalty0
  (2):\penalty0 489--504, 2001.

\bibitem[Bernard et~al.(2001)Bernard, B{\'e}lair, and
  Mackey]{bernard2001sufficient}
Samuel Bernard, Jacques B{\'e}lair, and Michael~C Mackey.
\newblock Sufficient conditions for stability of linear differential equations
  with distributed delay.
\newblock \emph{Discrete \& Continuous Dynamical Systems-B}, 1\penalty0
  (2):\penalty0 233, 2001.

\bibitem[Braverman and Zhukovskiy(2012)]{braverman2012absolute}
Elena Braverman and Sergey Zhukovskiy.
\newblock Absolute and delay-dependent stability of equations with a
  distributed delay.
\newblock \emph{Discrete \& Continuous Dynamical Systems-A}, 32\penalty0
  (6):\penalty0 2041, 2012.

\bibitem[Campbell and Jessop(2009)]{campbell2009approximating}
SA~Campbell and R~Jessop.
\newblock Approximating the stability region for a differential equation with a
  distributed delay.
\newblock \emph{Mathematical Modelling of Natural Phenomena}, 4\penalty0
  (2):\penalty0 1--27, 2009.

\bibitem[Cao(2019)]{cao2019bifurcations}
Yang Cao.
\newblock Bifurcations in an internet congestion control system with
  distributed delay.
\newblock \emph{Applied Mathematics and Computation}, 347:\penalty0 54--63,
  2019.

\bibitem[Cassidy(2020)]{cassidy2020distributed}
Tyler Cassidy.
\newblock Distributed delay differential equation representations of cyclic
  differential equations.
\newblock \emph{arXiv preprint arXiv:2007.03173}, 2020.

\bibitem[Doldo and Pender(2020)]{doldo2020multi}
Philip Doldo and Jamol Pender.
\newblock {Multi-Delay Differential Equations: A Taylor Expansion Approach}.
\newblock \emph{arXiv preprint arXiv:2012.05005}, 2020.

\bibitem[Jessop and Campbell(2010)]{jessop2010approximating}
R~Jessop and Sue~Ann Campbell.
\newblock Approximating the stability region of a neural network with a general
  distribution of delays.
\newblock \emph{Neural Networks}, 23\penalty0 (10):\penalty0 1187--1201, 2010.

\bibitem[Mor{\u{a}}rescu et~al.(2007)Mor{\u{a}}rescu, Niculescu, and
  Gu]{moruarescu2007stability}
Constantin-Irinel Mor{\u{a}}rescu, Silviu-Iulian Niculescu, and Keqin Gu.
\newblock Stability crossing curves of shifted gamma-distributed delay systems.
\newblock \emph{SIAM Journal on Applied Dynamical Systems}, 6\penalty0
  (2):\penalty0 475--493, 2007.

\bibitem[Novitzky and Pender(2020)]{novitzky2020queues}
Sophia Novitzky and Jamol Pender.
\newblock {Queues with Delayed Information: A Probabilistic Perspective}.
\newblock \emph{Cornell University, Ithaca NY}, 14853, 2020.

\bibitem[Novitzky et~al.(2019)Novitzky, Pender, Rand, and
  Wesson]{novitzky2019nonlinear}
Sophia Novitzky, Jamol Pender, Richard~H Rand, and Elizabeth Wesson.
\newblock Nonlinear dynamics in queueing theory: Determining the size of
  oscillations in queues with delay.
\newblock \emph{SIAM Journal on Applied Dynamical Systems}, 18\penalty0
  (1):\penalty0 279--311, 2019.

\bibitem[Novitzky et~al.(2020)Novitzky, Pender, Rand, and
  Wesson]{novitzky2020limiting}
Sophia Novitzky, Jamol Pender, Richard~H Rand, and Elizabeth Wesson.
\newblock Limiting the oscillations in queues with delayed information through
  a novel type of delay announcement.
\newblock \emph{Queueing Systems}, 95\penalty0 (3):\penalty0 281--330, 2020.

\bibitem[Pender et~al.(2017)Pender, Rand, and Wesson]{pender2017queues}
Jamol Pender, Richard~H Rand, and Elizabeth Wesson.
\newblock Queues with choice via delay differential equations.
\newblock \emph{International Journal of Bifurcation and Chaos}, 27\penalty0
  (04):\penalty0 1730016, 2017.

\bibitem[Tang(2004)]{tang2004oscillation}
XH~Tang.
\newblock Oscillation of first order delay differential equations with
  distributed delay.
\newblock \emph{Journal of mathematical analysis and applications},
  289\penalty0 (2):\penalty0 367--378, 2004.

\end{thebibliography}



\end{document}